\documentclass[12pt,leqno]{article}
\usepackage{amsfonts}
\usepackage{amsmath, amsthm, amsfonts, amssymb, color}
\usepackage{mathrsfs}
\usepackage{color}
\theoremstyle{definition}

\newcommand{\scr}[1]{\mathscr #1}
\definecolor{wco}{rgb}{0.5,0.2,0.3}

\numberwithin{equation}{section} \theoremstyle{remark}

\newcommand{\ua}{\uparrow}

\title{{\bf  Gradient Estimates and Applications for Neumann Semigroup  on Narrow Strip }\footnote{Supported in
 part by  NNSFC(11131003) and Laboratory of  Mathematics and Complex Systems} }
\author{
{\bf     Feng-Yu Wang  }\\
\footnotesize{ School of Mathematical Sciences,
Beijing Normal
University, Beijing 100875, China}\\
 \footnotesize{ Department of Mathematics,
Swansea University, Singleton Park, SA2 8PP, United Kingdom}\\
\footnotesize{  wangfy@bnu.edu.cn, F.-Y.Wang@swansea.ac.uk}}
\begin{document}
\allowdisplaybreaks
\def\R{\mathbb R}  \def\ff{\frac} \def\ss{\sqrt} \def\B{\mathbf
B}
\def\N{\mathbb N} \def\kk{\kappa} \def\m{{\bf m}}
\def\ee{\varepsilon}\def\ddd{D^*}
\def\dd{\delta} \def\DD{\Delta} \def\vv{\varepsilon} \def\rr{\rho}
\def\<{\langle} \def\>{\rangle} \def\GG{\Gamma} \def\gg{\gamma}
  \def\nn{\nabla} \def\pp{\partial} \def\E{\mathbb E}
\def\d{\text{\rm{d}}} \def\bb{\beta} \def\aa{\alpha} \def\D{\scr D}
  \def\si{\sigma} \def\ess{\text{\rm{ess}}}
\def\beg{\begin} \def\beq{\begin{equation}}  \def\F{\scr F}
\def\Ric{\text{\rm{Ric}}} \def\Hess{\text{\rm{Hess}}}
\def\e{\text{\rm{e}}} \def\ua{\underline a} \def\OO{\Omega}  \def\oo{\omega}
 \def\tt{\tilde} \def\Ric{\text{\rm{Ric}}}
\def\cut{\text{\rm{cut}}} \def\P{\mathbb P} \def\ifn{I_n(f^{\bigotimes n})}
\def\C{\scr C}      \def\aaa{\mathbf{r}}     \def\r{r}
\def\gap{\text{\rm{gap}}} \def\prr{\pi_{{\bf m},\varrho}}  \def\r{\mathbf r}
\def\Z{\mathbb Z} \def\vrr{\varrho} \def\ll{\lambda}
\def\L{\scr L}\def\Tt{\tt} \def\TT{\tt}\def\II{\mathbb I}
\def\i{{\rm in}}\def\Sect{{\rm Sect}}  \def\H{\mathbb H}
\def\M{\scr M}\def\Q{\mathbb Q} \def\texto{\text{o}} \def\LL{\Lambda}
\def\Rank{{\rm Rank}} \def\B{\scr B} \def\i{{\rm i}} \def\HR{\hat{\R}^d}
\def\to{\rightarrow}\def\l{\ell}\def\iint{\int}
\def\EE{\scr E}\def\no{\nonumber}
\def\A{\scr A}\def\V{\mathbb V}
\def\BB{\scr B}\def\Ent{{\rm Ent}}\def\x{\mathbf x}\def\y{\mathbf y} \def\z{\mathbf z} \def\hd{\mathbf d}

\maketitle

\begin{abstract} By using local and global versions of Bismut type derivative formulas, gradient estimates  are derived for   the Neumann semigroup on a narrow strip. Applications to functional/cost inequalities and heat kernel estimates are presented.  Since the narrow strip we consider is non-convex with zero injectivity radius,   and does not satisfy the volume doubling condition, existing results   in the literature  do not apply.
\end{abstract} \noindent
 \vskip 2cm

\section{Introduction}

Let $\phi_1,\phi_2\in C^2(\R^d)$ with
$\phi_1<\phi_2$ and $\lim_{|x|\to\infty} \{\phi_2(x)-\phi_1(x)\}=0.$ We   investigate gradient estimates and applications for the Neumann semigroup on the   strip
$$D:=\big\{(x,y)\in\R^{d+1}: \phi_1(x)\le y\le \phi_2(x)\big\}.$$
As the condition $\lim_{|x|\to\infty} \{\phi_2(x)-\phi_1(x)\}=0$ means that the strip is extremely narrow at infinity, we call $D$ a narrow strip.  This feature leads to essential difficulties in the study of the Neumann semigroup:\beg{enumerate} \item[(a)]  The domain  is non-convex with injectivity zero, so that existing results on gradient estimates and applications derived in \cite{W05,W10a,WPC,W10b,JW} using Li-Yau's maximum principle   and probabilistic arguments do not apply. \item[(b)] The domain does not satisfy the volume doubling condition, so that the argument for heat kernel estimates  developed by Grigoy'an (see \cite{GT,BCS}  and   references therein)  using the doubling condition  does not work. \end{enumerate} As far as we know, the study of gradient and heat kernel estimates for the Neumann semigroup on   a narrow strip remains new.

\

Let $L= \DD+Z$ for some $C^1$-vector field $Z$ on $\R^d$. We consider the Neumann semigroup $P_t$ generated by $L$ on the narrow strip $D$. Throughout the paper, we  use $\DD$ and $\nn$ to denote the Laplacian and the gradient operators on the underlying Euclidean space. The main tools of our study are  local/global   derivative formulas addressed in Section 4.    To apply these  formulas, we need the following conditions on $\phi_i (i=1,2)$ and $Z$.
\beg{enumerate} \item[(i)] $\phi_1,\phi_2\in C^2(\R^d),\ \phi_1<\phi_2,$ \
 $\lim\limits_{|x|\to\infty} \{\phi_2(x)-\phi_1(x)\}=0,\  \liminf_{|x|\to\infty} \<\nn \phi_1,\nn \phi_2\> >-1, $  and $ \<\nn \phi_1,\nn\phi_2\>(x)\le   ( |\nn\phi_1|^2\land|\nn\phi_2|^2)(x)$  for large $|x|>0$.
  \item[(ii)] $\sup\{\<\nn_v Z(x,y),v\>:\ v\in\R^{d+1}, |v|\le 1, (x,y)\in D\}<\infty.$
\item[(iii)] For $i=1,2$,
$$\sup\bigg\{\ff{(-1)^i\Hess_{\phi_i}(a,a)}{\phi_2-\phi_1}(x):\ x,a\in\R^d, |a|= 1\bigg\}<\infty.$$
\item[(iv)] For $i=1,2,$
$$\sup_{(x,y)\in D} \Big\{(\phi_i(x)-y)\{\DD\phi_i(x)+\<(\nn\phi_i(x),-1), Z(x,y)\>\} +|\nn\phi_i|^2(x)\Big\}<\infty.$$
\item[(v)] $\limsup\limits_{|x|\to\infty}\sup_{y\in [\phi_1(x), \phi_2(x)]}\dfrac {-L(\phi_2-\phi_1)}{\phi_2-\phi_1}(x,y)<\infty,\ \ \limsup\limits_{|x|\to\infty} |\nn\log (\phi_2-\phi_1)(x)|<\infty.$
\end{enumerate}
In the first condition of (v),  and also in the sequel, a function $\phi$ on $\R^d$ is naturally extended to    $\R^{d+1}$ by setting $\phi(x,y):=\phi(x), (x,y)\in\R^{d+1}.$

\

Under  these conditions,   the reflecting diffusion process generated by $L$ on $D$ is non-explosive. More precisely, consider the following stochastic differential equation with reflection:
\beq\label{E0} \d (X_t,Y_t) = \ss 2\, \d B_t +Z(X_t,Y_t)\d t + N(X_t,Y_t)\,\d l_t,\end{equation}
where $B_t$ is the $(d+1)$-dimensional Brownian motion, $N$ is the unit inward normal vector field of $\pp D$, and $l_t$ is the local time of the solution $(X_t,Y_t)$ on $\pp D$. Under the above conditions, for any initial data $(x,y)\in D$, the equation has a unique solution $\{(X_t,Y_t)(x,y)\}_{t\ge 0}$ which is non-explosive (see Proposition \ref{P2.1} below). Then the  Neumann semigroup generated by $L$ is  formulated as
$$P_t f(x,y)= \E f((X_t,Y_t)(x,y)),\ \ (x,y)\in D, t\ge 0, f\in \B_b(D).$$

\beg{thm}\label{T1.1} Assume {\rm (i)-(v)}. For any initial data $(X_0,Y_0)\in D$, the equation $\eqref{E0}$ has a unique solution which is non-explosive. Moreover, there exists a constant $c>0$ such that the associated Neumann semigroup $P_t$ satisfies the following gradient estimates. \beg{enumerate}
\item[$(1)$] For any $p>1$,
$$|\nn P_t f|\le (P_t|\nn f|^p)^{\ff 1 p} \exp\Big[c+\ff{cpt}{p-1} \Big],\ \ t\ge 0, f\in C_b^1(D).$$
\item[$(2)$] For any $p\in (1,2], t>0$ and $f\in \B_b(D),$
$$|\nn P_t f|^2\le  \ff{cp}{2(p-1)^2(1-\exp[-cpt/(p-1)])}(P_t |f|^p)^{\ff 2 p} .$$\end{enumerate}\end{thm}

Next, we present some applications of Theorem \ref{T1.1}. Let $\rr_D$ be the intrinsic distance on $D$, i.e. for any $\x,\y\in D$,
$$\rr_D(\x,\y):= \inf\bigg\{\int_0^1 |\gg'(s)|\d s: \gg\in C^1([0,1];D),\gg(0)=\x,\gg(1)=\y\bigg\}.$$ Moreover, for any probability measures $\mu$ and $\nu$ on $D$,
$$W_2^{\rr_D}(\mu,\nu):= \inf_{\pi\in \scr C(\mu,\nu)} \bigg(\int_{D\times D} \rr^2_D\d\pi\bigg)^{\ff 1 2}$$ is   the corresponding $L^2$-Wasserstein distance between $\mu$ and $\nu$, where $\scr C(\mu,\nu)$ is the set of all couplings of $\mu$ and $\nu$.
The following assertions  are more or less standard consequences of the gradient estimates in Theorem \ref{T1.1}.

\beg{cor}\label{C1.2} Assume {\rm (i)-(v)}. There exists a constant $c>0$ such that the following assertions hold. \beg{enumerate} \item[$(1)$] For any $t>0,$ the following Poincar\'e inequality holds:
$$P_t f^2 \le (P_t f)^2 + \ff{\e^{c}(\e^{ct}-1)}{c}P_t |\nn f|^2,\ \ f\in C_b^1(D).$$
\item[$(2)$] For any   $t>0$, the following log-Harnack inequality holds:
$$P_t (\log f)(\x) \le \log P_t f(\y) + \ff{c\e^{c}\rr_D(\x,\y)^2}{1-\e^{-ct}},\ \ \x,\y\in D, 0<f\in \B_b(D).$$
\item[$(3)$] For any measure $\mu$  which is  equivalent to the Lebesgue measure on $D$, the density $p_t(\x,\y)$ of $P_t$ with respect to $\mu$ satisfies the following entropy inequality:
$$\int_D p_t(\x,\z)\log \ff{p_t(\x,\z)}{p_t(\y,\z)}\,\mu(\d \z) \le  \ff{c\e^{c}\rr_D(\x,\y)^2}{1-\e^{-ct}},\ \ \x,\y\in D, t>0.$$
\item[$(4)$] The invariant probability measure  $\mu$ of $P_t$ is unique, and if it exists then the adjoint operator $P_t^*$ of $P_t$ on $L^2(\mu)$ satisfies the following entropy-cost inequality:
$$\int_D (P_t^*f)\log P_t^*f\d\mu\le  \ff{c\e^{c} }{1-\e^{-ct}}W_2^{\rr_D}(f\mu,\mu),\ \ t>0, f\ge 0,\mu(f)=1.$$
\item[$(5)$] Let $\mu$ be the invariant probability measure of $P_t$. Then the density $p_t(\x,\y)$ of $P_t$ with respect to $\mu$ satisfies
$$\int_D p_t(\x,\z)p_t(\y,\z)\mu(\d \z)\ge \exp\bigg[-\ff{\rr_D(\x,\y)^2  c\e^{c} }{1-\e^{-ct}}\bigg],\ \ t>0, \x,\y\in D.$$
If   $P_t$ is symmetric in $L^2(\mu)$, then
$$p_t(\x,\y)\ge \exp\bigg[-\ff{\rr_D(\x,\y)^2  c\e^{c} }{1-\e^{-ct/2}}\bigg],\ \ t>0, \x,\y\in D.$$\end{enumerate}\end{cor}

To illustrate the above results, we consider the following example where $\phi_i (i=1,2)$ are functions of $|x|$ for large $|x|$.

\paragraph{Example 1.1.} Let $\phi_i(x)= \ll_i \varphi(|x|)\ (i=1,2)$ for large $|x|$, where   $\ll_1<\ll_2$ with
 $\ll_1\le 0\le\ll_2$ are two constants, and $\varphi\in C_b^2([0,\infty))$ with $\varphi>0,   \varphi(r)\downarrow 0$ as $r\uparrow\infty$, and
  $$\limsup_{r\to\infty}\ff{|\varphi''(r)|+|\varphi'(r)|}{\varphi(r)}<\infty. $$ Moreover, let $Z= (Z_1,Z_2)\in C^1(\R^{d+1}; \R^d\times\R^1)$ satisfy (ii) and
$$\limsup_{|x|\to\infty}\sup_{y\in [\ll_1\varphi(x),\ll_2\varphi(x)]}\,    \Big(\ff{|\varphi'(|x|)|\cdot |Z_1(x,y)|}{\varphi(|x|)} + \varphi(|x|) |Z_2(x,y)|\Big) <\infty.$$
 Then it is easy to see that conditions (i)-(v) also hold, so that Theorem \ref{T1.1} and Corollary \ref{C1.2} apply.

 Typical choices of $\varphi(r)$ for large $r$ meeting the above requirements include $\varphi(r)= \e^{-\ll r^\dd}$ for some $\ll>0$ and $\dd\in (0,1]$, $\varphi(r)= r^{-\dd}$ for some $\dd>0$, and $\varphi(r)= \log^{-\dd}(\e+r)$ for some $\dd>0.$

\

The remainder of the paper is organized as follows. In Section 2,  we present some preparations, which include  the non-explosion of the reflecting diffusion process,  exponential estimates  on the local time, and
a prior gradient estimate on $P_t$. In Section 3,  we prove  Theorem \ref{T1.1} and Corollary \ref{C1.2}. Finally, in Section 4,  we introduce the local/global derivative formulas of the Neumann semigroup, which are used in Sections 2-3 as  fundamental tools.

\section{Preparations}

The main tool in our study of  gradient estimates   is the following   derivative formula (see Theorem \ref{A1} below):
\beq\label{BB} \nn P_t f=\E\big\{Q_t^* \nn f(X_t,Y_t)\big\} =\ff 1{\ss 2} \E\bigg\{f(X_t,Y_t) \int_0^t h_s' Q_s^*\d B_s\bigg\},\end{equation} where $h\in C^1([0,t])$ is such that $h_0=0$ and $h_t=1$, $Q_s$ is an adapted   process on $\R^d\otimes\R^d$ satisfying
\beq\label{BB'} \|Q_s\|\le \e^{\int_0^s K(X_r,Y_r)\d r+\int_0^s \si(X_r,Y_r)\d l_r},\end{equation} and $-\si$ is a lower bound of the second fundamental form of the boundary $\pp D$. So, to apply this formula, we need to calculate the second fundamental form,  and to estimate the exponential moment of the local time.
Moreover, to ensure the validity of the above derivative formula, we also need to prove the non-explosion of the reflecting diffusion process generated by $L$, and to verify the boundedness of $\nn P_t f$ on $[0,t]\times D$ for a reasonable class of functions $f$. These will be done in the following three subsections respectively.

\subsection{The second fundamental form}

Let $\pp_i=\{(x,\phi_i(x)): x\in\R^d\}$. We have $\pp D=\pp_1\cup\pp_2.$ Let $N$ be the unit inward normal vector field on $\pp D$. Then
\beq\label{N}N(x,\phi_i(x))= \ff {(-1)^i (\nn\phi_i(x),-1)}{\ss{1+|\nn \phi_i(x)|^2}},\ \ x\in\R^d, i=1,2.\end{equation} Recall that the second fundamental form of $\pp D$ is the following symmetric two-tensor defined on $T\pp D$, the tangent space of $\pp D$:
 $$\II(u,v):= - \<\nn_u N, v\>=-\<\nn_v N,u\>,\ \ \ u,v\in T\pp D.$$ We say that the second fundamental form is bounded below by a function $-\si$ on $\pp D$ and denote $\II\ge -\si$, if
$$\II(v,v)\ge -\si(\z) |v|^2,\ \ \ \z\in\pp D, v\in T_\z\pp D.$$ Below, we  calculate the lower bound of the second fundamental form.

For  any unit tangent vector $v$ of $\pp D$ at point $(x,\phi_i(x))\in\pp_i$,  there exists $a\in\R^d$ with $|a|=1$ such that
$$v= \ff{(a,\nn_a\phi_i(x))} {\ss{1+|\nn_a\phi_i(x)|^2}}.$$ Combining this with \eqref{N}, we obtain
$$\II(v,v)= -\<\nn_v N, v\>= -\ff{(-1)^{i}\Hess_{\phi_i}(a,a)}{\ss{1+|\nn\phi_i|^2} (1+|\nn_a\phi_i|^2)}(x).$$ Therefore, letting
\beg{equation}\label{SI}\beg{split} &\si_i(x)=  \sup_{a\in\R^d,|a|=1} \ff{(-1)^{i}\Hess_{\phi_i}(a,a)}{\ss{1+|\nn\phi_i|^2} (1+|\nn_a\phi_i|^2)}(x),\ \ x\in\R^d,i=1,2,\\
&\si(x,y)= \si_1(x) 1_{\{y=\phi_1(x)\}}+ \si_2(x) 1_{\{y=\phi_2(x)\}},\ \ (x,y)\in \pp D,\end{split}\end{equation}  we obtain
$ \II \ge -\si.$

\subsection{Non-explosion and exponential estimates on   $l_t$}

To investigate the non-explosion, we introduce the following Lyapunov function:
$$W_0(x,y)=W_0(x)= \ff 1{(\phi_2-\phi_1)(x)},\ \  (x,y)\in D.$$  By \eqref{N} and (i), there exists $r_0>0$ such that
\beq\label{NK} N W_0(x,\phi_i(x))=\ff{\<\nn \phi_2,\nn \phi_1\>-|\nn\phi_i|^2}{(\phi_2-\phi_1)^2\ss{1+|\nn\phi_i|^2}}(x)\le 0,\ \ i=1,2, W_0(x) \ge r_0.\end{equation} Thus, $NW_0\le 0$ holds on $\pp D\cap\{W_0\ge r_0\}.$ We modify $W_0$ such that this boundary condition holds on the whole boundary $\pp D$. Take $\bb\in C^\infty([0,\infty))$ with $\bb'\ge 0, \bb|_{[0,r_0]}=r_0,$ and $\bb(r)=r$ for $r\ge r_0+1.$ Then \eqref{NK} implies
\beq\label{NK'}    N W|_{\pp D}\le 0,\ \ W:=\bb\circ W_0.\end{equation}
Moreover, by (i), $W$ is a compact function on $D$, i.e. $\{\x\in D: W(\x)\le r\}$ is compact for any $r> 0.$ Define
\beq\label{tau} \tau_n=\inf\Big\{t\ge 0: \ W(X_t,Y_t)=W(X_t)\ge n\Big\},\ \ n\ge 1.\end{equation} Then the life time of the process can be formulated as
$$\tau_\infty =\lim_{n\to\infty}\tau_n.$$
\beg{lem}\label{P2.1} Assume {\rm (i)-(v)}. \beg{enumerate} \item [$(1)$] For any initial data $(X_0,Y_0)\in D$, the unique solution to the equation $\eqref{E0}$ is non-explosive.
\item[$(2)$]  For any $R>0$, there exists a constant $c>0$ such that for any initial data $(X_0,Y_0)\in D$,
$$\E \e^{\ll\int_0^t 1_{\{W\le R\}}(X_s)\d l_s}\le \e^{\ll + c\ll(1+\ll)t},\ \ \ll,t\ge 0.$$
\item[$(3)$] There exists a constant $c>0$ such that for any initial data $(X_0,Y_0)\in D$,
$$\E\e^{\ll \int_0^t \si(X_s,Y_s)\d l_s}\le \e^{c\ll+ c\ll(1+\ll)t},\ \ \ll, t\ge 0.$$\end{enumerate} \end{lem}

\beg{proof} (1) It is easy to see from (v) and the construction of $W$  that
$L W\le CW$ holds for some constant $C>0$.  Then by \eqref{NK'} and  It\^o's formula, we obtain
$$\d W(X_t,Y_t)\le \d M_t +CW(X_t,Y_t)\d t $$ for some local martingale $M_t$. This implies
$$\E W(X_{t\land\tau_n}, Y_{t\land\tau_n}) \le W(X_0,Y_0)\e^{Ct},\ \ t\ge 0.$$ Since $W>0$ and $W(X_{\tau_n\land t}, Y_{\tau_n\land t})=n$ provided $\tau_n\le t$, it follows that
$$\P(\tau_n\le t)\le \ff {\E W(X_{t\land\tau_n},Y_{t\land\tau_n})} n \le \ff{\e^{Ct}W(X_0,Y_0)}n.$$  Therefore, $\P(\tau_\infty\le t)=0$ holds for any $t\ge 0$, i.e. the process is non-explosive.

(2) Let   $\rr_\pp =\inf_{\y\in \pp D} |\cdot-\y|$ be the distance function to the boundary $\pp D$. Since $D$ is a $C^2$-domain, $\rr_\pp$ is $C^2$-smooth in a neighborhood of $\pp D$. So, for any $R>0$, there exists $\vv\in (0,1)$ such that
$\rr_\pp\in C^2\big(D\cap \{W\le R+1\}\cap \{\rr_\pp\le \vv\}\big).$ Let $\aa,\bb\in C^\infty([0,\infty))$ such that $\aa(r)=r$ for $r\le \ff{\vv}2,\ \aa|_{[\vv,\infty)}=\vv;$ and $\bb|_{[0,R]}=1,\ \bb|_{[R+1,\infty)}=0,\ \bb'\le 0.$ Then
$$\tt\rr_\pp:= (\bb\circ W) (\aa\circ\rr_\pp)\in C_0^2(D),\ \ 0\le \tt\rr_\pp \le 1.$$ Moreover, since on $\pp D$ we have $\rr_\pp=0$ and $ N\rr_\pp=1$, it follows from \eqref{NK'}, $\bb'\le 0, \aa'(0)=1$ and $\bb\circ W\ge 1_{\{W\le R\}}$ that
$$N\tt\rr_\pp = (\bb'\circ W) (\aa\circ\rr_\pp)NW +(\bb\circ W) \aa'(0)N\rr_\pp \ge 1_{ \{W\le R\}}$$ holds on $\pp D$. Thus, by It\^o's formula we obtain
\beg{equation*}\beg{split} \d \tt\rr_\pp(X_t,Y_t) &=  \ss 2\, \<\nn\tt\rr_\pp (X_t,Y_t), \d B_t\> + L \tt\rr_\pp(X_t,Y_t)\d t + N \tt\rr_\pp(X_t,Y_t)\d l_t\\
&\ge \ss 2\, \<\nn\tt\rr_\pp (X_t,Y_t), \d B_t\>- \|L\tt\rr_\pp\|_\infty\d t + 1_{\{W(X_t)\le R\}}\d l_t.\end{split}\end{equation*} Therefore,
$$\E \e^{\ll \int_0^t 1_{\{W(X_s)\le R\}}\d l_s} \le \e^{\ll +\ll \|L\tt\rr_\pp\|_\infty t} \E \e^{\ll \ss 2 \int_0^t \<\nn\tt\rr_\pp (X_s,Y_s), \d B_s\>} \le \e^{\ll + \ll \|L\tt\rr_\pp\|_\infty t+ \ll^2 \|\nn\tt\rr_\pp\|_\infty^2 t}.$$

(3) Let $g_i(x,y)= (\phi_i(x)-y)^2, (x,y)\in D.$ We have $(N g_{i})(x,\phi_{i}(x))=0$ (i.e. $N g_i|_{\pp_i}=0$) and
\beq\label{LL}Ng_i(x,\phi_{3-i}(x))=-\ff{2(\phi_2-\phi_1)(1+\<\nn\phi_1,\nn\phi_2\>)}{\ss{1+|\nn\phi_{3-i}|^2}}(x)=:-\tt\si_{3-i}(x),\ \ x\in\R^d.\end{equation} By (i), $1+\<\nn\phi_1,\nn\phi_2\>(x)\ge \theta_0$ holds for some constant $\theta_0>0$ and large enough $|x|>0$. Then it follows from \eqref{SI}, \eqref{LL} and (iii) that $\si_i\le \theta \tt\si_i$ holds on $\{W\ge R\}\cap \pp D$ for some constants $\theta, R>0.$ Since $\si_i$ is bounded on the compact set $\pp D\cap \{W\le R\},$ we conclude that
\beq\label{L*} \si_i\le \theta\tt\si_i +c_1 1_{\{W\le R\}}\end{equation} holds on $\pp D$ for some constant $c_1>0.$
   Moreover, by (iv) we have
\beg{equation*}\beg{split} Lg_{3-i}(x,y)&= 2-2y\DD\phi_{3-i}(x)+\DD \phi_{3-i}^2(x) + 2  (\phi_{3-i}(x)-y)\<(\nn \phi_{3-i}(x),-1), Z(x,y)\>\\
&\le K_{3-i},\ \ (x,y)\in D\end{split}\end{equation*} for some constant $K_{3-i}>0.$  Combining this with  $Ng_{3-i}|_{\pp_{3-i}}=0$ and \eqref{LL}, and using   It\^o's formula, we obtain
\beq\label{MT} \beg{split} & \d g_{3-i}(X_t,Y_t)\le \d M_t+ K_{3-i}\d t -\tt\si_i(X_t) \d l_t^i,\\
&M_t := 2 \ss 2\int_0^t (\phi_{3-i}(X_s)-Y_s)\<(\nn\phi_{3-i}(X_s),-1),\d B_s\>,\end{split}\end{equation}
where $l_t^i$ is the local time of $(X_t,Y_t)$ on $\pp_i$. Due to  \eqref{L*}, \eqref{MT} and    $(\phi_i(x)-y)^2\le\dd^2$   on $D$, we arrive at
\beg{equation*}\beg{split} \ll \int_0^t\si_i(X_s)\d l_s^i&\le  \ll\theta \int_0^{t} \tt\si_i(X_s)\d l_s^i+ \ll c_1 \int_0^t 1_{\{W(X_s)\le R\}}\d l_s^i\\
 &\le  \ll\theta  \big\{\dd^2+ K_{3-i}t\big\}  +\ll M_t+ \ll c_1 \int_0^t 1_{\{W(X_s)\le R\}}\d l_s^i.\end{split}\end{equation*}
Noting that
\beg{equation*}\beg{split}\<M\>_t &= 8 \ll^2\theta^2 \int_0^t |\phi_{3-i}(X_s)-Y_s|^2(1+|\nn\phi_{3-i}(X_s)|^2)\d s \\
&\le 8\ll^2\theta^2\dd^2(1+\|\nn\phi_{3-i}\|_\infty^2)t=:c_2\ll^2t,\end{split}\end{equation*}
this together with (2) implies
 \beg{equation*}\beg{split} \E\e^{\ll\int_0^t\si_i(X_s)\d l_s^i}&\le    \e^{\ll\theta(\dd^2+K_{3-i}t)+c_2\ll^2 t}\E^{M_t- \<M\>_t+\ll c_1 \int_0^t 1_{\{W(X_s)\le R\}}\d l_s^i} \\
 &= \e^{\ll\theta(\dd^2+K_{3-i}t)+c_2\ll^2t }\big(\E\e^{2M_t-2\<M\>_t}\big)^{\ff 1 2} \big(\E \e^{2\ll c_1 \int_0^t 1_{\{W(X_s)\le R\}}\d l_s^i}\big)^{\ff 1 2}\\
 &\le \e^{c\ll + c\ll(1+\ll)t}\end{split}\end{equation*} for some constant $c>0.$
 Therefore, we prove  (3) by noting that
 \beg{equation*}\beg{split} &\E\e^{\ll \int_0^t \si(X_s,Y_s)\d l_s}= \E\e^{\ll \int_0^t \si_1(X_s)\d l_s^1+\ll \int_0^t \si_2(X_s)\d l_s^2}\\
 &\le \Big(\E  \e^{2\ll\int_0^t\si_1(X_s)\d l_s^1}\Big)^{\ff 1 2} \Big( \E\e^{2\ll\int_0^t\si_2(X_s)\d l_s^2}\Big)^{\ff 1 2}.\end{split}\end{equation*}
\end{proof}

\subsection{A prior gradient estimate  on $P_t$}

In this subsection, we    prove the boundedness of   $\nn P_\cdot f$ on $[0,t]\times D$ for a nice reference function $f$ such that the derivative formula  \eqref{BB} is valid  according to Theorem \ref{A1}.  To this end, we use  the  local derivative formula presented in Theorem \ref{A2} below. The key point to apply this  formula lies in the construction of the   control process $h_s$, which is non-trivial due to the stopping time $\tau_B^x$. In  \cite[Section 4]{TW}, this control process was constructed by using a time change induced by the distance function to the boundary. However, in the present case  the distance of a point $\x\in D$ to the boundary vanishes as $|\x|\to\infty.$ So,   the construction  from \cite[Section 4]{TW} does  not imply  the desired boundedness of $\nn P_\cdot f$ on $[0,t]\times D$.    Our trick to fix this point is to use the Lyapunov function $W$ instead of the distance function to $\pp D$,  where $W$ is in \eqref{NK'}.

Let
$$g_n(x)=\cos\ff{\pi W(x)}{2n},\ \ x\in\R^d, n\ge 1.$$ For fixed $X_0\in\R^d$,  we consider     $n>W(X_0)+1+r_0,$ where $r_0>0$ is in \eqref{NK}. Define
$$T(t)= \int_0^tg_n(X_{s\land\tau_n})^{-2}  \d s,\ \ t\ge 0, $$ where $\tau_n$ is in \eqref{tau}. Then $T\in C([0,\tau_n);[1,\infty))$ is strictly increasing with $T(t)\ge t$, and $T(t)=\infty$ holds for $t>\tau_n$.   Let
$$\tau(t)= \inf\{s\ge 0: T(s)\ge t\},\ \ t\ge 0.$$ We have    $\tau(t)\le t,$ and $T\circ \tau(t)=t$ provided $\tau(t)<\tau_n$.

\beg{lem}\label{L2.1} Assume {\rm (i)-(v)}. Let $X_0\in \R^d$ and $n> W(X_0)+ r_0+1.$   Then $\tau(t)<\tau_n$ holds for all $t>0$. Moreover,   for any $m\ge 1$, there exists a constant $c>0$ independent of $n$  such that
\beq\label{B0} \E g_n(X_{\tau(t)})^{-m} \le g_n(X_0)^{-m} \e^{ct},\ \ t\ge 0.\end{equation} \end{lem}

\beg{proof} Let $\zeta_l=\inf\{t\ge 0: g_n(X_{\tau(t)})\le \ff 1 l\},\ l\ge 1.$ By the definitions of $g_n$ and $\tau_n$, we have
\beq\label{ZT} \zeta_\infty:=\lim_{l\to\infty} \zeta_l =\inf\{t\ge 0: \tau(t)\ge\tau_n\}.\end{equation} Moreover, by \eqref{NK'} we have $Ng_n^{-r}|_{\pp D}\le 0$ for any $r>0.$   So, by  It\^o's formula and the fact that
$$\d \tau(t)= g_n(X_{\tau(t)})^2\d t,\ \ t\le\zeta_l,$$ we obtain
\beq\label{I1} \beg{split} \d g_n(X_{\tau(t)})^{-m}&= \d M_t +(L g_n^{-m})(X_{\tau(t)},Y_{\tau(t)})\d\tau(t) + Ng_n^{-m}(X_{\tau(t)},Y_{\tau(t)})\d l_t\\
&=\d M_t + g_n(X_{\tau(t)})^2 (Lg_n^{-m})(X_{\tau(t)},Y_{\tau(t)})\d t,\ \ t\le \zeta_l\end{split}\end{equation} for some martingale $M_t$.
Since $W=W_0$ for $W_0\ge r_0+1$, by (v), there exists a constant $C>0$ independent of $n$ such that
\beg{equation*}\beg{split} &-(L g_n) (x,y)  = \bigg(\ff{\pi\sin\ff{\pi}{2n(\phi_2-\phi_1)}}{2n} \big(L(\phi_2-\phi_1)^{-1}\big)  +\ff{\pi^2g_n(x)}{4n^2}\big|\nn(\phi_2-\phi_1)^{-1}\big|^2\bigg)(x,y)\\
&= \bigg(\ff{\pi\sin\ff{\pi}{2n(\phi_2-\phi_1)}}{2n(\phi_2-\phi_1)} \Big( 2|\nn\log(\phi_2-\phi_1)|^2 -\ff{L (\phi_2-\phi_1)}{\phi_2-\phi_1}\Big) +\ff{\pi^2g_n|\nn\log(\phi_2-\phi_1)|^2}{4n^2 (\phi_2-\phi_1)^2}\bigg)(x,y) \\
&\le C,\ \  \text{if}\ (x,y)\in D,  1+r_0\le W_0(x)\le   n,\end{split}\end{equation*}and
$$|\nn g_n|^2(x) \le \ff{\pi^2|\nn\log(\phi_2-\phi_1)|^2}{4n^2(\phi_2-\phi_1)^2}(x)\le C,\ \ (x,y)\in D,  1+r_0\le W_0(x)\le   n.$$   Moreover, it is easy to see that $\{|Lg_n|+|\nn g_n|\}_{n\ge 1}$ are uniformly bounded on the compact set $D\cap \{W_0\le r_0+1\}$, we conclude that
\beg{equation*}\beg{split} &g_n(X_{\tau(t)})^2 (Lg_n^{-m})(X_{\tau(t)},Y_{\tau(t)})) \\
&= -m g_n(X_{\tau(t)})^{1-m} (Lg_n)(X_{\tau(t)},Y_{\tau(t)})) + m(m+1)(g_n^{-m}|\nn g_n|^2)(X_{\tau(t)})\\
 &\le c g_n(X_{\tau(t)})^{-m},\ \ t\le \tau_n\end{split}\end{equation*} holds for some constant $c>0$ independent of $n$. Combining this with \eqref{I1}, we obtain
\beq \label{BB}\E g_n(X_{\tau(t\land\zeta_l)})^{-m}\le g(X_0)^{-m}\e^{ct},\ \ t\ge 0, l\ge 1.\end{equation} Thus,
 $$\P(\zeta_l\le t)\le \ff{\E g_n(X_{\tau(t\land \zeta_l)})^{-m}}{l^m} \le \ff{\e^{ct}}{g_n(X_0)^{m} l^m}.$$ Letting $l\to\infty$ we obtain
 $\P(\zeta_\infty\le t)=0$ for all $t>0$, so that by \eqref{ZT}, $\P(\tau(t)<\tau_n)=1$ holds for all $t\ge 0$.  Finally, \eqref{B0} follows from \eqref{BB} by letting $l\to\infty$.
\end{proof}

Combining Lemma \ref{L2.1} with  Theorem \ref{A2} below, we can prove the following gradient estimate on $P_t$.

\beg{lem}\label{P2.2} Assume {\rm (i)-(v)}. There exists a constant $C>0$ such that
\beq\label{B22} |\nn P_t f|^2 \le \ff{C}{t\land 1} \big\{P_t f^2-(P_tf)^2\big\},\ \ t>0, f\in \B_b(D).\end{equation} Consequently, for any $f\in C_0^2(D)$ satisfying the Neumann boundary condition, $\nn P_\cdot f$ is bounded on $[0,t]\times D$ for all $t>0.$ \end{lem}

\beg{proof} We first observe that it suffices to prove \eqref{B22}. Indeed, if \eqref{B22} holds, then  for $f\in C_0^2(D)$ satisfying the Neumann boundary condition,
$$P_t f^2 -(P_t f)^2 = f^2 +\int_0^t (P_s Lf^2)\d s -\bigg(f+\int_0^t (P_s Lf)\d s \bigg)^2\le C(t+t^2),\ \ t\ge 0$$ holds for some constant $C>0.$ Combining this with \eqref{B2}, we conclude that $\nn P_\cdot f$ is bounded on $[0,t]\times D$ for all $t>0.$

Next,  by the semigroup property and Jensen's inequality, we only need to prove \eqref{B22} for $t\in (0,1].$ Moreover, by an approximation argument, for the proof of \eqref{B22} we may and do assume that  $f\in C_0^2(D)$ satisfying the Neumann boundary condition. Finally, for fixed $t>0$ and $\x_0\in D$, by using $f-P_tf(\x_0)$ to replace $f$,  to prove \eqref{B22} at point $\x_0$ we may assume further that $P_t f(\x_0)=0.$

Now, let $t\in (0,1], \x_0\in D$, and $f\in C_0^2(D)$ satisfy the Neumann boundary condition with $P_t f(\x_0)=0.$   For   $n> 1+r_0+W(\x_0),$ let
$$h_s= \ff 1 t \int_0^s g_n(X_r)^{-2}1_{\{r<\tau(t)\}}\d r,\ \ s\ge 0,$$ where $(X_t,Y_t)$ solves the equation \eqref{E0} with $(X_0,Y_0)=\x_0.$ Then $h_0=0$. Moreover, if $s\ge\tau(t),$ then
$$h_s=h_{\tau(t)}=\ff 1 t \int_0^{\tau(t)}g_n(X_r)^{-2}\d r= \ff{T\circ\tau(t)}t=1,$$where the last step follows since  $\tau(t)<\tau_n$ according to   Lemma \ref{L2.1}, so that $T\circ\tau(t)=t$ by the definitions of $T$ and $\tau$. By (ii), $\II\ge -\si$, and   $\tau(t)\le t\land\tau_n$, we can apply Theorem \ref{A2} below for $D_0=\{(x,y)\in D: W(x)\le n\}$ and
$$\|Q_s\|\le \e^{K(s\land\tau_n) +\int_0^{s\land\tau_n} \si(X_r,Y_r)\d l_r},\ \ s\ge 0.$$ By the Cauchy-Schwartz inequality, the definition of $h_s$,  Lemmas \ref{P2.1}-\ref{L2.1},  and $\tau(t)\le t$,   we obtain
\beg{equation*}\beg{split} &\E\bigg(\int_0^t |h'(s)|^2\|Q_s\|^2\d s\bigg)^2 \le \ff 1 {t^3} \E\int_0^{\tau(t)} g_n(X_s)^{-8} \e^{4Ks+4\int_0^{s\land\tau_n}\si(X_r)\d l_r}\d s\\
&\le\ff{\e^{4Kt}}{t^3}\bigg(\E\int_0^{\tau(t)} g_n (X_s)^{-16}\d s\bigg)^{\ff 1 2} \bigg(\E\int_0^{\tau(t)} \e^{8\int_0^s \si(X_r)\d l_r} \bigg)^{\ff 1 2}\\
&\le \ff{\e^{c_1(t+1)}}{t^{\ff 5 2}}\bigg(\E\int_0^{\tau(t)} g_n(X_s)^{-14}\d T(s)\bigg)^{\ff 1 2} \\
&\le \ff{\e^{c_1(t+1)}}{t^{\ff 5 2}} \bigg(\E\int_0^t g_n(X_{\tau(s)})^{-14}\d s\bigg)^{\ff 1 2} \le\ff{\e^{c_2(t+1)}}{t^2 g_n(\x_0)^7}\end{split}\end{equation*} for some constants $c_1,c_2>0$ independent of $n$.  Therefore, by Theorem \ref{A2} below,
$$|\nn P_tf|^2(\x_0)\le \ff 1 2(P_tf^2(\x_0))\E\int_0^t |h'(s)|^2 \|Q_s\|^2\d s \le \ff{C}{tg_n(\x_0)^{\ff 7 2}}P_t f^2(\x_0),\ \ t\in (0,1] $$ holds for some constant $C>0$ and all $n> 1+r_0+W(\x_0)$.  Then the proof is finished by letting $n\to\infty$.
\end{proof}

\section{Proofs of Theorem \ref{T1.1} and Corollary  \ref{C1.2}}

\beg{proof}[Proof of Theorem \ref{T1.1}] By an approximation argument, we may and do assume that    $f\in C_0^2(D)$ satisfying the Neumann boundary condition. In this case,  $\nn P_\cdot f$ is bounded on $[0,t]\times D$ according to   Lemma \ref{P2.2}.

(a) By (ii), $\II\ge -\si$,  and Theorem \ref{A1} below, we have
$$|\nn P_t f|\le \E\big\{|\nn f|(X_t,Y_t) \e^{Kt+\int_0^t\si(X_s,Y_s)\d l_s}\big\}\le \e^{Kt} (P_t|\nn f|^p)^{\ff 1 p} \big(\E\e^{\ff p{p-1} \int_0^t\si(X_s,Y_s)\d l_s}\big)^{\ff{p-1}p}.$$ Then the gradient estimate in (1) follows from Lemma \ref{P2.1}.

(b) Let $p\in (1,2]$ and $v\in\R^{d+1}$ with $|v|=1.$  By (ii),$\<\nn_v Z,v\>\le K|v|^2$ holds for some constant $K>0$ and all $v\in \R^{d+1}$. Combining this with  $\II\ge -\si$ and  using Theorem \ref{A1} below, for any $h\in C^1([0,t])$ with $h_0=0$ and $h_t=1$, we have
\beq\label{B1} |\nn_v P_t f|^2\le \ff 1 { 2} (P_t |f|^p)^{\ff 2 p} (\E |M_t|^q)^{\ff 2 q},  \ q:=\ff p{p-1}\ge 2,\end{equation}where
\beq\label{B2}M_t:= \int_0^t \<v,h'_sQ_s^*\d B_s\> \end{equation} for some adapted process $Q_s$ satisfying
\beq\label{B3} \|Q_s\|\le \e^{Ks+\int_0^s \si(X_r,Y_r)\d l_r},\ \ s\ge 0.\end{equation}
Noting that $\d|M_t|^2= 2 M_t \d M_t+ \d\<M\>_t$, for any $\vv>0$ we have
\beg{equation*}\beg{split} \d (M_t^2+\vv)^{\ff q 2} &= \d N_t + \ff q 2 (M_t^2+\vv)^{\ff q 2 -1} \d\<M\>_t + \ff q 4 \Big(\ff q 2 -1\Big) (M_t^2+\vv)^{\ff q 2 -2} 4M_t^2 \d\<M\>_t\\
&\le \d N_t + \ff{q(q-1)} 2 (M_t^2+\vv)^{\ff q 2 -1} \|h_t' Q_t\|^2 \d t, \end{split} \end{equation*}
where $\d N_T:= q M_t(M_t^2+\vv)^{\ff q 2 -1} \d M_t$ is a martingale due to Lemma \ref{P2.1} and \eqref{B2}. Moreover, letting   $\eta_t= \E (M_t^2+\vv)^{\ff q 2}$, and combining this with Lemma \ref{P2.1}, we obtain
\beg{equation*}\beg{split} \eta_t'&\le \ff {q(q-1)} 2\E \big\{(M_t^2+\vv)^{\ff q 2 -1} \|h_t' Q_t\|^2\big\}\le \ff {q(q-1)}2 |h_t'|^2 \eta_t^{\ff{q-2}q} (\E\|Q_t\|^q)^{\ff 2 q}\\
&\le \ff {q(q-1)}2 |h_t'|^2 \eta_t^{\ff{q-2}q}  \e^{2Kt} \big(\E\e^{q\int_0^t \si(X_s,Y_s)\d l_s}\big)^{\ff 2 q}\\
 &\le \ff {q(q-1)}2 |h_t'|^2 \eta_t^{\ff{q-2}q}   \e^{c+cqt}\end{split}\end{equation*} for some constant $c>0$. Therefore,
$$\eta_t\le \bigg(\vv^{\ff 2 q} + (q-1)\e^{c} \int_0^t  |h_s'|^2 \e^{cqs}\d s\bigg)^{\ff q 2}.
$$ Letting $\vv\to 0$ and taking
\beq\label{H} h_0=0,\ \ h_s'= \ff{\exp[-cqs]}{\int_0^t \exp[-cqs]\d s},\ \ s\in [0,t],\end{equation} we arrive at
$$\E |M_t|^q \le   \bigg(\ff{cq(q-1)\e^{c}}{1-\exp[-cqt]}\bigg)^{\ff q 2}.$$ Substituting this into \eqref{B1}, we prove (2).
\end{proof}

\beg{proof}[Proof of Corollary \ref{C1.2}] By Theorem \ref{T1.1}(1) with $p=2$, we have
\beq\label{BG} |\nn P_t f|^2\le (P_t|\nn f|^2)  \e^{c+ct},\ \ t\ge 0, f\in C_b^1(D) \end{equation}for some constant $c>0.$
By an approximation argument, in (1) and (2) we may and do assume that $f\in C_b^2(D)$ satisfying the Neumann boundary condition, which is constant outside a bounded set.

(a) The desired Poincar\'e inequality follows from \eqref{BG} and the following simple calculations due to Bakry-Emery (cf. \cite{Bakry}):
\beg{equation*}\beg{split} & P_tf^2-(P_tf)^2 =-\int_0^t \ff{\d}{\d s} P_{t-s} (P_s f)^2\d s  =\int_0^t P_{t-s}|\nn P_s f|^2 \d s\\
 &\le P_t |\nn f |^2 \int_0^t\e^{c+cs} \d s= \ff{\e^{c}(\e^{ct}-1)}{c}P_t |\nn f|^2.\end{split}\end{equation*}

(b) The proof of (2) can be modified from that of Theorem 2.1 in \cite{RW10}. More precisely, let $\x,\y\in D$ and $t>0$ be fixed. By the definition of $\rr_D(\x,\y)$, for any $\vv>0$, there exists a $C^1$-curve   $\gg: [0,1]\to D$    such that
\beq\label{CC} \gg(0)=\x, \ \gg(1)=\y,\ |\gg'|\le\rr_D(\x,\y)+\vv.\end{equation}  Let
$$h_0=0,\ \ h_s'= \ff{c \e^{-c s}}{1-\e^{-c t}},\ \ s\in [0,t]. $$ If $f\in C_b^2(D)$ is positive, constant outside a compact set,  and satisfies the Neumann boundary condition,  then
\beg{equation*}\beg{split} &P_t\log f(x)-\log P_t f(y)= \int_0^t \ff{\d}{\d s} (P_s\log P_{t-s}f)(\gg(h_s))\d s\\
&= \int_0^t \big\{h_s' \<\gg'(s),\nn P_s\log P_{t-s} f\>-P_s |\nn \log P_{t-s}f|^2\big\}(\gg(h_s))\d s\\
&\le\int_0^t \big\{|h_s'|(\rr_D(\x,\y)+\vv)|\nn P_s\log P_{t-s} f| -\e^{-c-cs}|\nn P_s\log P_{t-s}f|^2\big\}(\gg(h_s))\d s\\
&\le \ff{(\rr_D(\x,\y)+\vv)^2}4 \int_0^t |h_s'|^2\e^{c+cs}\d s= \ff{c\e^{c}(\rr_D(\x,\y)+\vv)^2}{1-\e^{-ct}}.\end{split}\end{equation*} Then the desired log-Harnack inequality follows by letting $\vv\to 0$.

(c) Applying \cite[Lemma 3.1(4)-(5)]{WY} for $P=P_t$ and $\Phi(s)= \e^s$, the desired entropy inequality in (3) as well as the heat kernel estimate in (5) follow from (2). Moreover, the entropy-cost inequality in (4) follows from (2) and \cite[Corollary 1.2(3)]{RW10}.
\end{proof}

\section{Derivative formulas for $P_t$}

In this section, we introduce derivative formulas of $P_t$ on a $C^2$-domain $D$ in $\R^\hd$ for $d\ge 2$, where $P_t$ is the Neumann semigroup generated by $L:=\DD+Z$ on $D$ for some $C^1$-vector field $Z$. Let $K\in C(D)$ such that
\beq\label{CD} \<\nn_v Z(x),v\>\le K(x) |v|^2,\ \ x\in D, v\in\R^\hd.\end{equation} Consider the following stochastic  differential equations:
\beq\label{APE} \d X_t= \ss 2\,\d B_t +Z(X_t)\d t +N(X_t)\d l_t,\end{equation} where $B_t$ is the $\hd$-dimensional Brownian motion, $N$ is the inward unit normal vector field of $\pp D$, and $l_t$ is the local time of the solution on $\pp D$. We assume that for any $x\in D$, the solution $(X_t^x, l_t^x)_{t\ge 0}$ to this equation starting at $x$ is non-explosive. Then   the associated Neumann semigroup is formulated as
$$P_t f(x)= \E f(X_t^x),\ \ t\ge 0, x\in D, f\in \B_b(D).$$
Moreover, let $\si\in C(\pp D)$ such that
\beq\label{S} \II (v,v):= -\<\nn_v N(x), v\> \ge -\si(x) |v|^2,\ \ x\in\pp D, v\in T_x\pp D.\end{equation}
To state the derivative formulas, we introduce the class $$\scr C_N(D):=\big\{f\in C^2(D), Nf|_{\pp D}=0, Lf\in \B_b(D)\big\}.$$ By the Kolmogorov equations, we have
(see e.g. \cite[Theorem 3.1.3]{Wbook})
\beq\label{KL} \ff{\d }{\d t} P_t f= P_t Lf= LP_t f,\ \ t\ge 0, f\in \scr C_N(D).\end{equation}

The following global derivative formula is essentially taken from \cite[Theorem 3.2.1]{Wbook}. This type of derivative formula was proved by  Bismut \cite{Bismut} and Elworthy-Li \cite{EL}   on manifolds without boundary.

\beg{thm}\label{A1} Let $t>0$ and $x\in D$ be fixed. If
\beq\label{DP}  \sup_{s\in [0,t]}\E\e^{\int_0^s K(X_r^x)\d r +\int_0^s \si(X_r^x)\d l_r^x}<\infty,\end{equation}
then there exists an adapted $\R^\hd\otimes\R^\hd$-valued process $(Q_s)_{s\in [0,t]}$ with
\beq\label{DE} Q_0=I,\ \ \|Q_s\|\le \e^{\int_0^s K(X_r^x)\d r +\int_0^s \si(X_r^x)\d l_r^x},\ \ s\in [0,t],\end{equation} such that for any $f\in \scr C_N(D)$ with bounded $\nn P_\cdot f$ on $[0,t]\times D$,
\beq\label{BS} \nn P_t f(x)= \E\{Q_t^* \nn f(X_t^x)\} = \ff 1 {\ss 2}\E\bigg\{f(X_t^x) \int_0^t h_s'Q_s^*\d B_s\bigg\}\end{equation} holds for any
$h\in C^1([0,t])$ with $h_0=0$ and $h_t=1.$ \end{thm}

\beg{proof} The construction of $Q_s$ as well as the first equality in \eqref{BS} are essentially due to \cite{Hsu}. Once $Q_s$ is   constructed, the second equality in \eqref{BS} can be proved as in   \cite{T}.

(a) For any $z\in \pp D$, let $P_\pp^z$ be the projection onto the tangent space $T_z\pp D$ of $\pp D$ at point $z$. We have
$$P_\pp^z a= a-\<a,N(z)\>N(z),\ \ a\in \R^\hd.$$ Next, let $\II_\pp^z\in \R^\hd\otimes\R^\hd$ such that
$$\<\II_\pp^z a,b\>= \II_\pp(P_\pp^z a, P_\pp^z b),\ \ a,b\in\R^\hd.$$ Moreover, for any $n\ge 1$,  let $(Q_s^{(n)})_{s\ge 0}$ solve the following equation on $\R^\hd\otimes\R^\hd$:
\beq\label{AP}\beg{split}  \d Q_s^{(n)} = &(\nn Z(X_s^x))Q_s^{(n)}\d s\\
 &-\II^{X_s^x}Q_s^{(n)}\d l_s^x -(n+\si(X_s^x)^+) ((Q_s^{(n)})^*N(X_s^x))\otimes N(X_s^x)\,\d l_s^x,\ \ Q_0^{(n)}=I,\end{split}\end{equation}  where for any $v_1,v_2\in\R^\hd$, $v_1\otimes v_2\in\R^\hd\otimes\R^\hd$ is defined by
$$(v_1\otimes v_2)a= \<v_1,a\>v_2,\ \ a\in \R^\hd.$$ Then for any $a\in \R^\hd$, it follows from \eqref{CD}, \eqref{S} and \eqref{AP} that
\beg{equation}\label{AP0}\beg{split} \d|Q_s^{(n)}a|^2 &= 2\big\<\nn_{Q_s^{(n)}a}Z(X_s^x), Q_s^{(n)}a\big\> \d s- 2 \II(P_\pp^{X_s^x}Q_s^{(n)}a,P_\pp^{X_s^x}Q_s^{(n)}a) \d l_s^x \\
&\qquad - 2(n+ \si(X_s^x)^+)\<Q_s^{(n)}a, N(X_s^x)\>^2\d l_s^x \\
&\le  2|Q_s^{(n)}a|^2 \big\{\si(X_s^x)\d l_s^x +K(X_s^x)\d s\big\} -2n \<Q_s^{(n)}a, N(X_s^x)\>^2\d l_s^x.\end{split}\end{equation}
In particular,
\beq\label{AP1} \|Q_s^{(n)}\|^2 \le\e^{2 \int_0^s K(X_r^x)\d r +2 \int_0^s \si(X_r^x)\d l_r^x}<\infty,\ \ s\ge 0, n\ge 1.\end{equation} By \eqref{DP} and \eqref{AP1}, we obtain
$$\sup_{s\in [0,t]} \E \sup_{n\ge 1} \|Q_s^{(n)}\|<\infty.$$ So, that the sequence $\{Q_\cdot^{(n)}\}_{n\ge 1}$ is uniformly integrable in $L^1([0,t]\times\OO\to \R^\hd\otimes\R^\hd; \d s\times \P)$, and $\{Q_t^{(n)}\}_{n\ge 1}$ is uniformly integrable in $L^1(\OO\to\R^\hd\otimes\R^\hd;\P)$. Therefore,  there exists a subsequence $n_k\uparrow\infty$ and a progressively measurable process $(Q_s)_{s\in [0,t]}$ satisfying \eqref{DE} such that for any bounded measurable function $\xi: [0,t]\to \R^\hd$ and any bounded $\hd$-dimensional random variable $\eta$,
\beq\label{AP3} \lim_{k\to\infty} \E \int_0^t \big\<(Q_s^{(n_k)}-Q_s)a,\xi_s\big\>\d s  = \lim_{k\to\infty} \E \big\<Q_t^{(n_k)}-Q_t)a,\eta\big\>=0,\ \ a\in\R^\hd.\end{equation}
 Moreover,  for any $m\ge 1$ and
$\tau_m:= \inf\{t\ge 0: |X_t^x|\ge m\}$, it follows from \eqref{AP0} and \eqref{AP1} that
\beq\label{AP2} \beg{split} &\lim_{n\to\infty} \E \int_0^{t\land\tau_m} \big<(Q_s^{(n)})a, N(X_s^x)\big\>^2\d l_s^x\\
&\le \lim_{n\to\infty} \ff {|a|^2} {n}\bigg(1+ \E \int_0^{t\land\tau_m} \|Q_s^{(n)}\|^2\{|K|(X_s^x)\d s+|\si|(X_s^x)\d l_s^x\}\bigg) =0,\end{split}\end{equation} where in the last step we have used  \eqref{AP1}, the boundedness of $K$ and $\si$ on the compact set $D_m:=\{z\in D: |z|\le m\}$ (where we take $\si=0$ outside $\pp D$), and
\beq\label{EXP} \E \e^{\ll l_{t\land\tau_m}^x}<\infty,\ \ x\in D, t\ge 0, \ll>0, m\in\N\end{equation} according to the proof of \cite[Theorem 6.1]{W10a}. In fact, as in Lemma \ref{P2.1}(2), we have the stronger conclusion that $\E\e^{\ll \int_0^t 1_{\{|X_s^x|\le m\}}\d l_s^x}<\infty.$

(b) By \eqref{APE}, \eqref{KL}, \eqref{AP} and   that $P_{t-s}f$ satisfies the Neumann boundary condition, It\^o's formula yields
\beq\label{AP4}\beg{split} & \d (\nn_{Q_s^{(n)}a} P_{t-s}f)(X_s^x) =  \ss 2\, \big\<\nn(\nn_{A_s^{(n)}a} P_{t-s} f)(X_s^x), \d B_s\big\>\\
&\quad +\big\{\big\<\nn(\nn_{Q_s^{(n)}a}P_{t-s}f),N\big\>(X_s^x)
 - \II(P_\pp^{X_s^x}Q_s^{(n)}a, \nn P_{t-s}f)(X_s^x)\big\}\d l_s^x.\end{split}\end{equation}
Since $P_{t-s}f$ satisfies the Neumann boundary condition,   for any $z\in\pp D$ and $v\in T_z\pp D$ we have
 $$0= \<v, \nn\<N,\nn P_{t-s}f\>(z)\> = \Hess_{P_{t-s}f}(v,N)(z)- \II(v, \nn P_{t-s}f(z)).$$ Then  whenever $X_s^x\in \pp D$,
 \beg{equation*}\beg{split} &\big\<\nn(\nn_{Q_s^{(n)}a}P_{t-s}f),N\big\>(X_s^x) = \Hess_{P_{t-s}f}(N, Q_s^{(n)}a)\\
 &= \Hess_{P_{t-s}f}(N,N)(X_s^x) \<Q_s^{(n)}a, N(X_s^x)\> +\II(P_\pp^{X_s^x}Q_s^{(n)}a, \nn P_{t-s}f)(X_s^x).\end{split}\end{equation*} Combining this with \eqref{AP4}, we obtain
 \beq\label{AP5} \beg{split} &\d (\nn_{Q_s^{(n)}a} P_{t-s}f)(X_s^x)\\
  &=  \ss 2\, \big\<\nn(\nn_{Q_s^{(n)}a} P_{t-s} f)(X_s^x), \d B_s\big\> + \Hess_{P_{t-s} f}(N,N)(X_s^x) \<Q_s^{(n)}a, N(X_s^x)\>\d l_s^x.\end{split}\end{equation}
Since $P_\cdot f\in C_b^2([0,t]\times D_m)$ due to the compactness of $D_m$,    there exists a constant $C_m>0$ such that
\beq\label{APn}\beg{split}& \limsup_{k\to\infty}\E \int_0^{t\land\tau_m} \Big|  \Hess_{P_{t-s} f}(N,N)(X_s^x) \<Q_s^{(n)}a, N(X_s^x)\>\Big|\d l_s^x\\
&\le C_m (\E l_{t\land\tau_m}^x)^{\ff 1 2}\limsup_{k\to\infty}\bigg(\E \int_0^{t\land\tau_m}  \<Q_s^{(n)}a, N(X_s^x)\>^2\d l_s^x\bigg)^{\ff 1 2}=0,\end{split}\end{equation} where the last step follows from \eqref{EXP} and \eqref{AP2}. Moreover, since $\nn P_\cdot f$ is bounded on $[0,t]\times D$, it follows from \eqref{DP}, \eqref{AP1} and \eqref{AP3} that
\beq\label{APm} \beg{split} &\limsup_{m\to\infty}\limsup_{k\to\infty} \Big|\E\<\nn P_{t-t\land \tau_m} f(X_{t\land\tau_m}^x), Q_{t\land\tau_m}^{(n_k)} a\>-\E\<\nn f(X_t^x),Q_ta\>\Big|\\
&\le    \limsup_{k\to\infty} \Big|\E\<\nn   f(X_{t}^x), (Q_{t}^{(n_k)}-Q_t) a\>\Big|\\
&\qquad
+2 \sup_{[0,t]\times D} |\nn P_\cdot f| \limsup_{m\to\infty} \E \Big\{1_{\{\tau_m<t\}}\e^{\int_0^t K(X_s^x)\d s+ \int_0^t \si(X_s^x)\d l_s^x}\Big\}=0.\end{split}\end{equation} Combining
  \eqref{AP5}, \eqref{APn} and \eqref{APm},    we arrive at
\beg{equation*}\beg{split} (\nn_a P_t f) (x)&=  \lim_{m\to\infty} \lim_{k\to\infty}  \E \<\nn P_{t-t\land\tau_m}f(X_{t\land \tau_m}^x), Q_{t\land\tau_m}^{(n_k)}a\>\\
&\qquad -
\lim_{m\to\infty} \lim_{k\to\infty} \E\int_0^{t\land\tau_m} \Hess_{P_{t-s} f}(N,N) \<Q_s^{(n_k)}a, N(X_s^x)\>\d l_s^x \\
 &= \E\<\nn f(X_t^x),Q_t a\>.\end{split}\end{equation*} So,  the first equality in \eqref{BS} holds.

 (c) By \eqref{KL} and It\^o's formula we have  $$\d P_{t-s} f(X_s^x) = \ss 2\, \<\nn P_{t-s} f(X_s^x),\d B_s\>.$$ Then
 $$f(X_t^x)= P_t f(x) +\ss 2 \int_0^t \<\nn P_{t-s} f(X_s^x),\d B_s\>,$$ so that by \eqref{AP3},
 \beg{equation}\label{AG1}\beg{split}  &\ff 1 {\ss 2} \E\bigg(f(X_t^x) \int_0^t h'_s \<Q_s a, \d B_s\>\bigg)
 = \E\int_0^t h_s' \<Q_s a,\nn P_{t-s} f(X_s^x)\>\d s\\
 &= \lim_{k\to\infty} \E \int_0^t h_s' \<Q_s^{(n_k)}a, \nn P_{t-s} f(X_s^x)\>\d s.\end{split}\end{equation}
   Similarly to \eqref{APm}, we have
 \beg{equation}\label{APg}\beg{split} &\lim_{k\to\infty} \E \int_0^t h_s' \<Q_s^{(n_k)}a, \nn P_{t-s} f(X_s^x)\>\d s\\
 &=\lim_{m\to\infty} \lim_{k\to\infty} \int_0^t h_s' \E \<Q_{s\land\tau_m}^{(n_k)}a, \nn P_{t-s\land\tau_m} f(X_{s\land \tau_m}^x)\>\d s.\end{split}\end{equation} Finally, by \eqref{AP5} and \eqref{APn},
\beg{equation*}\beg{split} &\lim_{k\to\infty} \E \<Q_{s\land\tau_m}^{(n_k)}a, \nn P_{t-s\land\tau_m} f(X_{s\land \tau_m}^x)\>\\
&= (\nn_a P_t f)(x) + \lim_{k\to\infty}
 \E\int_0^{s\land\tau_m} \Hess_{P_{t-r} f}(N,N) \<Q_r^{(n_k)}a, N(X_r^x)\>\d l_r^x\\
 &= (\nn_aP_t f)(x)\end{split}\end{equation*} holds uniformly in $s\in [0,t].$
  Combining this with \eqref{AG1} and \eqref{APg}, we obtain
 $$\ff 1 {\ss 2\,} \E\bigg(f(X_t^x) \int_0^t h'_s \<Q_s a, \d B_s\>\bigg)= \int_0^t h_s' (\nn_a P_t f)(x)\d s=\nn_a P_t f(x).$$
\end{proof}

To verify the boundedness of $\nn P_\cdot f$ on $[0,t]\times D$ required in Theorem \ref{A1}, one may use the following local version of derivative formula, which is essentially due to \cite[Lemma 3.2.2]{Wbook}. This type of derivative formula   goes back to \cite{T} for  manifolds without boundary. For a compact subset $B$ of $D$, we let
$$\tau_B^x=\inf\{t\ge 0: X_t^x\notin B\},\ \ x\in B.$$

\beg{thm}\label{A2} Let $x\in D$, and let $B$ be a compact subset of $D$ such that ${\rm dist}(x, D\setminus B)>0.$ Then there exists an adapted process $(Q_s)_{s\in [0,t]}$ on $\R^\hd\otimes\R^\hd$ with $Q_0=I$ and
\beq\label{C0} \|Q_s\|\le \e^{\int_0^{s\land\tau_B^x}K(X_r^x)\d r +\int_0^{s\land\tau_B^x} \si(X_r^x)\d l_r^x},\ \ s\in [0,t],\end{equation} such that for any adapted process $(h_s)_{s\in [0,t]}$ with $h_0=0$,   $h_s=1$ for $s\ge t\land \tau_B^x$, and
\beq\label{HH} \E \int_0^t |h_s'|^2\|Q_s\|^2\d s<\infty,\end{equation} there holds
\beq\label{BS2} \nn P_t f(x)= \ff 1 {\ss 2 } \E\bigg(f(X_t^x) \int_0^t h_s' Q_s^*\d B_s\bigg),\ \ f\in \scr C_N(D).\end{equation}\end{thm}

\beg{proof} Let $Q_s^{(n)}$ solve \eqref{AP}. By \eqref{AP1} and the exponential integrability of $l_{t\land\tau_B^x}$ due to \cite[Theorem 6.1]{W10a}, we have
$$\sup_{n\ge 1} \E\int_0^{t} \|Q_{s\land \tau_B^x}^{(n)}\|^2\d s <\infty.$$ Then there exists a subsequence $n_k\uparrow\infty$ and a progressively measurable process $(Q_s)_{s\in [0,t]}$ such that
\beq\label{DDE}\lim_{k\to\infty} \E\int_0^t (Q_{s\land\tau_B^x}-Q_s)\xi_s\d s =0,\ \ \xi_\cdot\in L^2([0,t]\times \OO;\d s\times\P).\end{equation} Next, by \eqref{APn} for $\tau_B^x$ in place of $\tau_m^x$, and using \eqref{AP5}, we have
$$\lim_{k\to\infty}\E \int_0^t (1-h_s) \d\<\nn P_{t-s}f(X_s),Q_s^{(n_k)}a\>=0.$$ Combining this with \eqref{DDE},
$h_s'=0$ for $s\ge t\land\tau_B^x$, and noting that
$$\d (P_{t-s}f)(X_s^x)= \ss 2\, \<\nn P_{t-s} f (X_s^x), \d B_s\>,\ \ s\in [0,t],$$ we obtain
\beg{equation*}\beg{split} &\ff 1 {\ss 2} \E\bigg(f(X_t^x) \int_0^t \<h_s' Q_sa, \d B_s\>\bigg) = \E\int_0^t \<h_s' Q_s a,\nn P_{t-s} f(X_s^x)\>\d s \\
&=\lim_{k\to\infty} \E \int_0^t (h_s-1)' \<Q_{s\land \tau_B^x}^{(n_k)}a,\nn P_{t-s} f(X_s^x)\>\d s \\
&=\lim_{k\to\infty} \E\bigg(\<\nn P_{t-s}f (X_s^x), Q_{s\land\tau_B^x}^{(n_k)}a\>(h_s-1) \big|_0^t-\int_0^t (1-h_s)  \d\<\nn P_{t-}f(X_s),Q_s^{(n_k)}a\>\bigg)\\
&=  \<\nn P_t f(x),a\>.\end{split}\end{equation*}  \end{proof}

\paragraph{\bf Acknowledgement.} The author would like to thank Professor Lixin Yan   for helpful communications.

\beg{thebibliography}{99}

\bibitem{Bakry} D. Bakry,  \emph{L$'$hypercontractivit\'e et son utilisation en
th\'eorie des semigroups,} in  $``$Lectures on Probability Theory'',  Lecture Notes in Math.  1581, Springer-Verlag, Berlin, pp. 1--114, 1994.

\bibitem{Bismut} J.-M. Bismut,  \emph{Large Deviations and the Malliavin Calculus,} Birkh\"auser, Boston, MA, 1984.

\bibitem{BCS} S. Boutayeb, T. Coulhon,  A. Sikora, \emph{A new approach to poitwise heat kernel upper bounds on doubling metric measure spaces,}  arXiv:1311.0367.

\bibitem{EL} K.-D. Elworthy, X.-M.  Li, \emph{Formulae for the
derivatives of heat semigroups,} J. Funct. Anal. 125(1984),
  252--286.

\bibitem{GT} A. Grigor'yan, T. Andras, \emph{Two-sided estimates of heat kernels on metric spaces,}  Ann. Probab. 40(2012), 1212--1284.

\bibitem{Hsu} E. P. Hsu,  \emph{ Multiplicative functional for the heat
equation on manifolds with boundary,}  Michigan Math. J.  20(2002),  351--367.

\bibitem{T} A. Thalmaier, \emph{On the differentiation of heat semigroups and
Poisson integrals,}  Stoch. Stoch. Rep. 61(1997), 297--321.

\bibitem{TW} A. Thalmaier, F.-Y. Wang,  \emph{Gradient estimates
for harmonic functions on regular domains in Riemannian
manifolds,} J. Funct. Anal. 155(1998), 109--124.

\bibitem{RW10} 	M. R\"ockner, F.-Y. Wang, \emph{Log-Harnack  Inequality for Stochastic differential equations in Hilbert spaces and its consequences,} Inf. Dimens. Anal. Quant. Probab.. Relat. Top. 13(2010), 27--37.

\bibitem{W05} F.-Y.
Wang,  \emph{ Gradient estimates and the first
Neumann eigenvalue on manifolds with boundary,}  Stoch. Proc. Appl.
115(2005), 1475--1486.

\bibitem{W10a} F.-Y. Wang, \emph{ Harnack inequalities on manifolds with boundary and applications,}   J. Math. Pures Appl. 94(2010), 304--321.

\bibitem{WPC}	F.-Y. Wang, \emph{Gradient and Harnack inequalities on noncompact manifolds with boundary,}  Pacific J. Math. 245(2010), 185--200.

\bibitem{W10b} F.-Y. Wang, \emph{ Semigroup properties for the second fundamental form,} Docum. Math. 15(2010), 543--559.

\bibitem{Wbook} F.-Y. Wang, \emph{Analysis for Diffusion Processes on Riemannian Manifolds,} World Scientific, Singapore, 2013.

\bibitem{WY} F.-Y. Wang, C. Yuan, \emph{Harnack inequalities for functional SDEs with multiplicative noise and applications,}  Stochastic Process. Appl. 121(2011), 2692--2710.

\bibitem{JW} J. Wang, \emph{Global heat kernel estimates,} Pacific J. Math.   178(1997),  377--398.

\end{thebibliography}

\end{document}